\documentclass{article}
\usepackage[utf8]{inputenc}
\usepackage{amsmath,amssymb,amsthm,amsfonts}
\usepackage{authblk}
\usepackage{pgf,tikz}
\usepackage{biblatex}
\usepackage[hyperfootnotes=false]{hyperref}
\usepackage{cleveref}
\usepackage{algorithm}
\usepackage{algpseudocode}
\usepackage[hyphenbreaks]{breakurl}
\algrenewcommand\algorithmicrequire{\textbf{Input:}}
\algrenewcommand\algorithmicensure{\textbf{Output:}}

\newcommand{\BB}{\textbf{B}}
\newcommand{\Char}[1]{\ensuremath{\characteristic(#1)}}
\newcommand{\multideg}[1]{\ensuremath{\deg(#1)}}
\newcommand{\Vector}[1]{\ensuremath{\boldsymbol{#1}}}
\newcommand{\wlexl}[1]{>_{#1}}
\newcommand{\lexl}{>_{\text{lex}}}
\newcommand{\LM}[1]{\ensuremath{\text{LM}(#1)}}
\newcommand{\LC}[1]{\ensuremath{\text{LC}(#1)}}
\newcommand{\LT}[1]{\ensuremath{\text{LT}(#1)}}
\newcommand*\quot[2]{{^{\textstyle #1}\big/_{\textstyle #2}}}

\DeclareMathOperator*{\characteristic}{char}

\newtheorem{theorem}{Theorem}
\numberwithin{theorem}{section}

\newtheorem{proposition}[theorem]{Proposition}

\theoremstyle{definition}
\newtheorem{definition}[theorem]{Definition}
\newtheorem{example}[theorem]{Example}
\newtheorem{remark}[theorem]{Remark}

\newcommand\blfootnote[1]{%
  \begingroup
  \renewcommand\thefootnote{}\footnote{#1}%
  \addtocounter{footnote}{-1}%
  \endgroup
}

\newcommand{\subjclass}[2]{%
\blfootnote{#1 \emph{Mathematics subject classification.} #2}}%

\title{Coding Reliability with {\sc Aclus}-- \\ Did I correctly characterize my observations?}
\author[1]{Marcus Weber}
\affil[1]{Zuse Institute Berlin (ZIB), Takustr. 7, D-14195 Berlin, Germany, \newline email: weber@zib.de}
\author[2]{O\u{g}uzhan Y\"ur\"uk}
\affil[2]{Chair of Discrete Mathematics/Geometry, Technische Universit{\"a}t Berlin, Stra{\ss}e des 17. Juni 136, 10623 Berlin, Germany, email: yuruk@math-tu.berlin.de}
\date{\today}

\begin{document}

\maketitle

\begin{abstract}
    Describing observations or objects in non-mathematical disciplines can often be accomplished by answering a list of questions. These questions can be formulated in such a way that the only possible answers always are ``yes'' or ``no''. This article is about automatically checking such given binary data sets for inconsistencies and about finding possible logical rules valid for the analyzed objects.
\end{abstract}

\subjclass{2010}{62H30, 62P25, 91C20, 13P10}

\section{Introduction}

This article presents a method (Algebraic clustering, {\sc Aclus}) to check for possible inconsistencies in a given binary data table. A binary data table represents an incidence relation between a set of observations and a list of characteristics. Such data structures are used by many disciplines to store various information such as the data collected in a social studies survey, or the data collected in a controlled experiment. There are several ways of cross-checking the correctness of coded observations \cite{Connor2020}: Colleagues can be consulted to check the assignments (intercoder reliability, intercoder consistency). The coder can also carry out the coding again at a later point in time and examine deviations (intracoder reliability). In this article, we propose a new approach, namely to check whether certain classifications ``destroy'' simple rules about the observed objects.  

The data stored in a binary table can be represented as a set of polynomials in a Boolean polynomial ring. Boolean rings and algebras are classically studied by George Boole in mid 1800s \cite{boole1847mathematical}, but they are quite fundamental objects in various research areas outside mathematics such as computer science \cite{Shannon:Synthesis1949,Braun:DigitalComputerDesign}, biology  \cite{VelizJarrahLaubenbacher:PolyAlgebraInBio:2010}, or data analysis \cite{Chikalov:LAD2013, Alexe:Coronary2003,schwieger2020classifier}. Due to their practicality in various applications, the Boolean polynomial rings and their algebraic aspects are intensively studied in mathematics. 

The method {\sc Aclus} which we introduce in this article relies on a well-established theory in the field of computational algebra, called the theory of Gröbner bases. Computing a Gröbner basis is usually not an easy task for a general input, but there are highly efficient computation methods for Boolean polynomials. A Gröbner basis of a Boolean ideal can be computed efficiently by using zero suppressed binary decision diagrams (ZDDs) which exploit the extra structure of Boolean polynomial rings and can handle a few hundred variables \cite{Brickenstein:PolyBori}. Throughout this work, we use the Gröbner basis implementation available in the computer algebra system OSCAR \cite{OSCAR}, which is based on the software {\sc Singular} \cite{DGPS}. For even larger scale problems, one possibility is to get help from modern SAT solving techniques \cite{Biere2009HandbookSAT}. The performance of modern SAT-solvers are constantly thanks to the International SAT Competition organized annually, see \cite{SAT_Website}. However, the scale of our computations in this work are relatively small, and we will not need such heavy tools.

{\sc Aclus} is based on an existing theory in mathematics, however, this article adds novel concepts to data analysis in two different ways: 
\begin{itemize}
    \item The presented framework allows for representing ``objects'' as well as ``properties'' in the same space, namely by Boolean polynomials. This is in the spirit of the ``Mathematics of Comparing Objects'' \cite{WbF_2022}. Objects are represented by polynomials which combine properties. To extract the properties from the polynomials, one has to compare objects (i.e., adding and multiplying the according polynomials) which again leads to new polynomials. 
    \item  From a statistical point of view (like in machine learning) ``rules'' about properties of objects become visible by finding repeating patterns or significant regularities in data. In this article we think of ``logical rules''. From this point of view, a ``rule'' becomes visible by objects (or observations) that do {\em not} exist, i.e.  by the ``holes'' in our binary data tables expressed by an ideal $\cal I$. 
\end{itemize}

\subsection{Guiding Example}\label{sec:intro}
The guiding example of our article is inspired by a decision tree in \cite{rotbuche}, and it is given in \Cref{tab:my_label}. Such tables are called binary data sets or Boolean data tables.
\begin{definition}
\label{Definition:BooleanDataTable}
A Boolean data table (binary data set) $M$ is a matrix whose entries are either 0 or 1. In general, we assume that $M$ has $m$ rows which correspond to the observations and $n$ columns which correspond to the characteristics.
\end{definition}

Our guiding example analyzes under which conditions the copper beech has a healthy growth. Every row of this table represents an experiment in which a certain soil and temperature condition is arranged. Afterwards it is checked whether the copper beech has a good growth. 

\begin{table}[t]
    \centering
    \begin{tabular}{|l||c|c|c|c|c||c|}
         \hline
         No. & wet (w)& acidic (a) & neutral (n) & cold (c) & very (v) & grow (g)  \\
         \hline
         \hline
         1 &0 & 0 & 0 & 0& 0& 1 \cr
         2 &0 & 0 & 0 & 1& 0& 1 \cr
         3 &0 & 0 & 0 & 1& 1& 1 \cr
         4 &0 & 0 & 1 & 0& 0& 0 \cr
         5 &0 & 0 & 1 & 1& 0& 1 \cr
         6 &0 & 0 & 1 & 1& 1& 1 \cr
         7 &0 & 1 & 0 & 0& 0& 1 \cr
         8 &0 & 1 & 0 & 1& 0& 1 \cr
         9 &0 & 1 & 0 & 1& 1& 0 \cr
         10& 1 & 0 & 0 & 0& 0& 1 \cr
         11& 1 & 0 & 0 & 1& 0& 0 \cr
         12&1 & 0 & 0 & 1& 1& 0 \cr
         13 &1 & 0 & 1 & 0& 0& 1 \cr
         14&1 & 0 & 1 & 1& 0& 0 \cr
         15&1 & 0 & 1 & 1& 1& 0 \cr
         16&1 & 1 & 0 & 0& 0& 1 \cr
         17&1 & 1 & 0 & 1& 0& 0 \cr
         18&1 & 1 & 0 & 1& 1& 0 \cr
         \hline
    \end{tabular}
    \caption{Binary data set. This table could come from a laboratory notebook in which observations of the growth of copper beeches were depicted. It is reported whether the soil is wet (w), has an acidic pH (a), or a neutral pH (n). It is also reported whether the environment is cold (c) or even very cold (v). Based on these experimental conditions it is checked whether the trees can grow well (g).}
    \label{tab:my_label}
\end{table}

What is special about line No. 4 in Table \ref{tab:my_label}? In this case the environment of the tree is not wet, the tree has a neutral soil, it is warm, and the tree does not grow well. At a first glance, it is remarkable that it is not cold and despite this circumstance, the tree does not grow well. Why is this special? In any other case, a not cold environment leads to a well growing tree. The algebraic approach in this article will reveal an even more precise statement: In every other trial (except No. 4) at least one of the two characteristics -- ``cold ($c$)'' or ``the tree is growing well ($g$)'' -- is assigned to the observation. This means: If observation No. 4 would not be there, then the condition ``it is not cold'' would be {\em sufficient but not necessary} for the good growth of the copper beech, i.e. then the rule would be: If not $c$ then $g$.

There are now two possibilities that can lead to this particular observation No. 4. On the one hand, it may really be a special feature of the copper beech that it does not always grow well if it is warm (because it is dry). On the other hand, however, this entry in the laboratory book could also simply be wrong and the corresponding experiment should perhaps be repeated.

\subsection{Mathematical Problem} 
{\em Given a binary data table whose entries are not random, but contain (hidden) logical rules which have led to the special combination of assigned and not-assigned characteristics of the observations. How can I check whether this table includes assignment errors contradicting the unknown rules? }

When analyzing a binary data table, checking how often a certain pattern occurs is a crucial step. Consider the pattern given in the line No. 5, which is ``001101''. If the same pattern occurs in many other lines in the data table, then it would make sense to trust this particular entry since it is not an exception in the data. However, this is not the criterion of correctness in this article. In fact, {\sc Aclus} will first collect the set of all pairwise different patterns in the table and then only work on these different patterns. From this point of view, two observations are not distinguishable if they share the same pattern. In Table \ref{tab:my_label} every pattern only occurs once. Thus, the step of clustering lines sharing the same pattern is obsolete in this special example. If a logical rule concerning the properties of the given objects applies, then the negation of this rule provides a not satisfiable statement.  For example, take the logical rule: If it rains, then the soil becomes wet. Its negation -- it rains and the soil stays dry -- provides a situation that will not be observed. The set of all patterns which can not be found among our observations is therefore the basis for finding logical rules about our reported object properties. In order to formulate this approach algebraically, we represent these logic and set operations in Boolean rings. 

A ring with unity $1\in \BB$ is called a \emph{Boolean ring} if every element $a\in \BB$ is idempotent, i.e., $a^2 = a$. It follows from the the idempotency assumption that any Boolean ring $\BB$ has characteristic 2, i.e. $a+a = 0$ for any $a\in \BB$. Moreover, for any $a,b \in \BB$ we have $ a+b = (a+b)^2 = a+a \cdot b+b \cdot a+b $, and since $\Char{\BB}=2$ it follows that $\BB$ is commutative.

Operations in Boolean rings can be seen as logical statements. Given any Boolean ring $(\BB, + , \cdot)$, one can define new set of Boolean operations $\lor, \land, \lnot$ (i.e., ``or'', ``and'', ``not'') where 
\begin{align*}
	a \lor b = a + b + a \cdot b, \quad a \land b = a \cdot b \quad \text{and} \quad \lnot a = a + 1, 
\end{align*}
for $a,b \in \BB$. The set $\BB$ together with these 3 new Boolean operations defines a Boolean algebra. Conversely, if $(\BB, \lor , \land,\lnot)$ is a Boolean algebra then we can recover the ring addition $a + b = \left( \lnot a \land b\right) \lor \left(a \land \lnot b  \right)$ and multiplication $a \cdot b = a \land b$ for any $a,b \in \BB$. 

%There is a natural way to define a partial order $\succcurlyeq$ on any Boolean ring $\BB$. Namely, we write $a \succcurlyeq b$ for $a,b \in \BB$ if and only if $ab = b$. In terms of Boolean operators over sets, this translates to $a \land b = b$ which is equivalent to saying that $b$ is contained in $a$. 
	
The operations of Boolean algebras are also understandable as operations on sets. Given an arbitrary set $S$, its power set $\mathcal{P}(S)$ has a natural Boolean algebra structure.
If two subsets $A,B \in \mathcal{P}(S)$, then their union and intersection correspond to binary Boolean operations $\lor$ and $\land$, respectively, and the complement $A^{c} := S \setminus A$ is the unary operation $(1+A)$ associated to the Boolean algebra $\mathcal{P}(S)$. The set $S$ does not necessarily have to be finite, but in the scope of our application this set will always be finite. 
The empty set and $S$ are distinguished elements of the Boolean algebra $\mathcal{P}(S)$, because they are the identity of the intersection and the union operations. The ring operations induced by the Boolean algebra structure is then as follows:
	$$ P + Q = (P \cap Q^c) \cup (P^c \cap Q) \text{ and } P \cdot Q = P \cap Q. $$ 
For more details about the Boolean rings and algebras, we refer the reader to \cite{Halmos:LecturesOnBooleanAlgebras}.	
%Instead of analyzing single patterns which do not occur in the reported table, we analyse all possible algebraic expressions in Boolean rings which provide subsets with  non-observed patterns. 
 
\begin{example}
The algebraic approach presented in this article works out for each row of a binary table in which way this row is ``special''. It creates an algebraic expression, which names the special feature of the corresponding pattern. To anticipate the presented method and the evaluations: In our example in line No. 4 it provides $cg+c+g+1$, which is equivalent to $\lnot(c\vee g)$. From the algebraic expression one can read which additional rule about the relationships between the columns would hold, if the corresponding pattern would be deleted from the table. In our example this rule would be $cg+c+g+1=0$, which is equivalent to $cg+c+g=1$, to $\lnot{c}\Rightarrow g$, or to $\lnot{g}\Rightarrow c$.

This Boolean algebraic approach will further turn into an algebraic clustering method which provides a way to not only find single, conspicuous rows in the table, but also how to combine some rows to get a conspicuous cluster of observations together with its peculiarity. In our guiding example, a combination of line No. 3 and line No. 6 forms one cluster. In every other row, ``very cold ($v$)'' and ``well-growing ($g$)'' do not appear together. Without these entries: $v\Rightarrow \lnot{g}$, equivalently $g\Rightarrow \lnot{v}$. 
\end{example}

\section{The Algebraic Framework}

In this section we present the algebraic framework behind our method {\sc Aclus}.  

Let $\BB$ be a Boolean ring, and $\BB[\Vector{X}] := \BB[X_1,\dots,X_n]$ denote the ring of polynomials over $\BB$. Every polynomial in $\BB[\Vector{X}]$ is a finite sum of monomials, and a \emph{monomial} in $\BB[\Vector{X}]$ is an expression of the form 
$$ f_i X_1^{\alpha_1} \dots X_n^{\alpha_n},$$ 
where $f_i \in \BB$ and $\Vector{\alpha} = (\alpha_1, \dots, \alpha_n) \in \mathbb{Z}_{\geq 0}^n$ are called as the coefficient and the exponent of this particular monomial.  $\BB[\Vector{X}]$ is not a Boolean ring, as the polynomial variable $X_i$'s are not necessarily idempotent. However, if the idempotency is enforced on the variables by considering the quotient ring 
	$$\BB(\Vector{X}) =  \quot{\BB[\Vector{X}]}{\left\langle X_1^2 - X_1, \dots, X_n^2-X_n \right\rangle}, $$
then the resulting quotient is a Boolean ring, and consequently every finitely generated ideal in $\BB(\Vector{X})$  is principal. An element from the quotient $\BB(\Vector{X})$ is called a Boolean polynomial, and it has an unique representative in $B[X_1, \dots , X_n]$ with degree at most 1 in each variable $X_i$. 

Let $\mathcal{O}$ be a set of observations, $\mathcal{P}$ be a set of characteristics associated to the elements in $\mathcal{O}$. Each observation in $\mathcal{O}$ either posses or lacks each characteristic in $\mathcal{P}$, and the incidence information between the observations and the characteristics can be represented in terms of a Boolean data table as in \Cref{Definition:BooleanDataTable}. Each characteristic in $\mathcal{P}$ is assigned a Boolean function (variable) $X_i$ from $\mathcal{O}$ to $\BB = \{0,1\}$. For a given observation $o \in \mathcal{O}$, the Boolean function $X_i$ is evaluated to 1 if and only if the observation $o$ admits the characteristic represented by $X_i$. 

Algebraically, the Boolean polynomial ring $\BB(X_1, \dots, X_n)$ contains every Boolean polynomial that can be written using the variables $X_1, \dots, X_n$. In terms of the Boolean data, $\BB(X_1, \dots, X_n)$ represents all possible descriptions that can be given using the characteristics represented by $X_1, \dots, X_n$. For example, $X_i \cdot X_j$ is a description which holds for any observation that admits both of the characteristics $X_i$ and $X_j$, and similarly, $X_i + X_j$ describes those observations that admit the characteristics $X_i$ or $X_j$ but not both. 

Given a Boolean polynomial $f(\Vector{X}) \in \BB(\Vector{X})$ and a vector $\Vector{b_o} \in \BB^n$, one can \emph{evaluate} the polynomial $f$ at $\Vector{b_o}$. Note that fixing a vector $\Vector{b_o} \in \BB^n$ corresponds to considering a particular observation in $\mathcal{O}$. In algebraic terms, this corresponds to substituting each variable $X_i$ with the $i$-th entry of $\Vector{b_o}$, and then computing the resulting Boolean algebraic expression.

\begin{example}
Consider the Table \ref{tab:my_label} again. Every row of this table can be mapped to a Boolean polynomial. Line No. 5 turns into $(w+1)(a+1)nc(v+1)(g+1)$. If we insert logical values ``0'' and ``1'' into this expression (in terms of Boolean ring elements), then $(w+1)(a+1)nc(v+1)(g+1)=1$ only for line No. 5. This algebraic expression is like a ``select statement'' for this special row. 
\end{example}

There are expressions in $\BB(\Vector{X})$ which are evaluated to ``0'' for all row patterns of Table \ref{tab:my_label}, e.g., $wan(c+1)vg=0$ for all lines in \Cref{tab:my_label}. If we have such an expression that is always evaluated to zero for all rows in the table, then every multiple of this expression (multiplication with Boolean polynomials) is also always zero. Similarly, this accounts for sums of such expressions, too.  In fact, this defines an ideal $\cal I$ of $\BB(\Vector{X})$, which we will refer as the \emph{Aclus ideal} associated to the given Boolean data table. In general, the Aclus ideal $\cal I$ includes all select-statements which would provide an empty set when applied to the data of the binary table. 
%If we understand that every logical rule about the relationship between the characteristics automatically defines a theoretical set of observations (a set of patterns) which are not to be found in ``reality'', then we can see that $\cal I$ provides all ``physical'' logical rules which apply to the reported characteristics.

\begin{example}
The Boolean polynomial $cv+v$ is element of the Aclus ideal $\cal I$ in case of Table \ref{tab:my_label}. The equation $cv+v=0$ is equivalent to $v=cv$. The logical rule is: If it is very cold, then it is cold and very cold. The expression $(g+1)w(c+1)$ is another element of this ideal, which means that whenever it is wet and not cold the copper beech grows, because $w(c+1)=1$ implies $g+1=0$ in order to cancel the algebraic expression. A generator of the Aclus ideal $\cal I$ can be found by computing the union of all rows of Table \ref{tab:my_label} and then by taking the complement of this expression. If we only consider the distinct patterns of Table \ref{tab:my_label}, then computing the union is like taking the sum of all related polynomial expressions. 
\end{example}

Given a Boolean polynomial ring $\BB(\Vector{X})$ and any ideal ${\cal I}$ of $\BB(\Vector{X})$, the set
	$$ V({\cal I}):=\left\{ \Vector{X}\in \BB^n \ \mid \ \forall f\in {\cal I}, \ f(\Vector{X}) = 0  \right\}$$
is called the variety of $\cal I$. Recall that an ideal is a collection of polynomials from $\BB(\Vector{X})$, and each of these polynomials represents a rule given in terms of the Boolean variables $X_1,\dots,X_n$. Then, $V({\cal I})$ is the set of objects in $\BB^n$ that do not satisfy any rule represented by the polynomials in $\cal I$.

Various fundamental results from classical algebraic geometry, such as \Cref{Theorem:BooleanNullstellensatz}, can be carried over to the Boolean setting. 
\begin{theorem}[Boolean Nullstellens\"{a}tze \cite{Sato:BooleanGroebner2011}]
	\label{Theorem:BooleanNullstellensatz}
	Let $\cal I$ be a finitely generated ideal in $\BB(\Vector{X})$, then $V({\cal I}) \neq \emptyset$ if and only if $\cal I$ contains a constant polynomial that is not equal to $0$. Moreover, given any $\cal I$ with $V({\cal I}) \neq \emptyset$ and $f \in \BB(\Vector{X})$, then 
		$$f\in {\cal I} \quad \Longleftrightarrow \quad f(\Vector{\tilde{X}}) = 0 \quad \text{for all} \quad \Vector{\tilde{X}} \in V({\cal I}).$$ 
\end{theorem}
\noindent In terms of Boolean tables, \Cref{Theorem:BooleanNullstellensatz} can be summarized by the following two statements: First, any rule generated from the polynomials in the ideal $\cal I$ is not satisfied by any element of $V({\cal I})$. Second, if a rule is not satisfied by any element of $V({\cal I})$, then this rule exists as a polynomial in $\cal I$. 

Given an ideal ${\cal I}$ of the Boolean polynomial ring $\BB(\Vector{X})$, a natural way to obtain a new polynomial ring is to consider the quotient ring $\quot{\BB(\Vector{X})}{\cal I}$. This quotient ring is also a Boolean ring, and its elements are not Boolean polynomials, but equivalence classes of Boolean polynomials in $\BB(\Vector{X})$. If we fix a suitable order on the monomials of $\BB(\Vector{X})$, then there is a way to pick a unique representative Boolean polynomial from each equivalence class. Therefore, one can associate the ring $\quot{\BB(\Vector{X})}{\cal I}$ with a new set of representative polynomials, which obviously depend on the fixed monomial order and this dependence will be carefully discussed in \Cref{Subsection:Grobner}. The elements of the ring $\quot{\BB(\Vector{X})}{\cal I}$ can be perceived as the reduced versions of the polynomials in $\BB(\Vector{X})$ such that each polynomial that lies in $\cal I$ is set to be equal to zero. For general polynomial rings, computing these representatives may be a challenging task. However, for Boolean polynomial rings, the theory of Gröbner bases allows us to do this efficiently. In \Cref{Subsection:Grobner}, we cover the basic notions about Boolean Gröbner basis which will be required to understand our approach.

\subsection{Simplifying Rules via Boolean Gröbner Basis}
\label{Subsection:Grobner}

Although in the ring of Boolean polynomials every finitely generated ideal is principal, there may be an alternative set of algebraic expressions which generate $\cal I$. 
%Finding a set of basic rules which apply for the given complete data set is then equal to finding a convenient basis of $\cal I$ consisting of expressions like $(v+1)c$ or $(g+1)(c+1)w$ for Table \ref{tab:my_label}.
One possible way to find a basis for the ideal $\cal I$ is to use the theory of Gröbner bases. Before we give a definition of Gröbner basis, we first have to introduce the notion of monomial orderings.

Recall that the set of monomials in $\BB({\Vector{X}})$ can be represented with the set $\mathbb{Z}^n_{\geq 0}$ by considering the exponent vector of each monomial in $\BB({\Vector{X}})$.
\begin{definition}
\label{Definition:Lexorder}
    A monomial ordering on $\BB(\Vector{X})$ is a relation $>$ on the set of exponent vectors $\mathbb{Z}^n_{\geq 0}$ which is a total and well ordering such that for any $\Vector{\alpha},\Vector{\beta}, \Vector{\gamma} \in \mathbb{Z}^n_{\geq 0}$ with $\Vector{\alpha}>\Vector{\beta}$, it holds that $\Vector{\alpha} + \Vector{\gamma}>\Vector{\beta} +\Vector{\gamma}$. 
\end{definition}
There are various ways to define a monomial order on $\BB({\Vector{X}})$, and here we present two standard examples. 
\begin{definition}[Lexicographic order]
Given two exponent vectors $$\Vector{\alpha} = (\alpha_1, \dots, \alpha_n), \Vector{\beta}= (\beta_1,\dots,\beta_n) \in \mathbb{Z}^n_{\geq 0},$$ we write $\Vector{\alpha} \lexl \Vector{\beta}$, or equivalently  $\Vector{X}^{\Vector{\alpha}} \lexl \Vector{X}^{\Vector{\beta}}$, if the leftmost nonzero entry of the vector $\Vector{\alpha}-\Vector{\beta} = (\alpha_1 - \beta_1 , \dots , \alpha_n-\beta_n)$ is positive. 
\end{definition}

\begin{example}
\label{Example:LexOrder}
Consider the monomials $X_1$, $X_2 X_3$, $X_2X_4$ and $X_2$ in the polynomial ring $\BB(X_1,X_2,X_3,X_4)$, whose exponent vectors are $(1,0,0,0)$,  $(0,1,1,0)$,  $(0,1,0,1)$ and $  (0,1,0,0)$, respectively. First, consider the lexicographic monomial order on $\BB(X_1,X_2,X_3,X_4)$: The difference of the first two exponent vectors is 
$$ (1,0,0,0) - (0,1,1,0) = (1,-1,-1,0), $$
and therefore $X_1 \lexl X_2X_3$. The differences of other exponent vectors similarly yield that $X_1 \lexl X_2 X_3 \lexl X_2X_4 \lexl X_2$. Furthermore, the lexicographic order on $\BB(X_1,X_2,X_3,X_4)$ induces the following order on the variables: 
\begin{equation}
    \label{Eqn:InducedOrderbyLex}
    X_1 \lexl X_2 \lexl X_3 \lexl X_4.
\end{equation}
\end{example}

The natural way to represent the exponent vector of a monomial is to write the exponent of $i$-th variable in the $i$-th entry of the vector. We make this assumption without exactly saying it in \Cref{Definition:Lexorder}, and this assumption induces the order of variables given in \eqref{Eqn:InducedOrderbyLex}. By changing in which order we write the variables into a vector, one can obtain the other possible ways to order variables. In fact, the variables are not necessarily given with a natural order in our applications. In \Cref{Algorithm:Weight}, we propose a methodology to obtain an variable order which makes sense from a data perspective. Also, it is possible to equip the variables with a preferred order by using the weighted lexicographic order, which we define next.

\begin{definition}[Weighted lexicographic order]
\label{Definition:WLexOrder}
Let $\Vector{w} = \left(w_1,\dots,w_n \right) \in \mathbb{Z}^n$ be a weight vector, and let $\lexl$ denote the lexicographic order on $\mathbb{Z}^n_{\geq 0}$. Then, the weighted lexicographic order given by $\Vector{w}$, which is denoted by $\wlexl{\Vector{w}}$, is defined as follows: for $\Vector{\alpha}=\left(\alpha_1,\dots,\alpha_n \right),\Vector{\beta}=\left(\beta_1,\dots,\beta_n \right)  \in \mathbb{Z}^n_{\geq 0}$, we write $\Vector{X}^{\alpha} \wlexl{\Vector{w}} \Vector{X}^{\Vector{\beta}}$ if and only if 
$$   \sum w_i\alpha_i > \sum w_i\beta_i \qquad \text{or} \qquad \sum w_i\alpha_i = \sum w_i\beta_i \ \ \text{and} \ \ \Vector{X}^{\alpha} \lexl \Vector{X}^{\Vector{\beta}}.$$
\end{definition}

\begin{example}
\label{Example:WLexOrder}
Consider the Boolean ring $\BB(X_1,X_2,X_3,X_4)$, the weight vector $\Vector{w} = (4,3,2,1)$, and let $\wlexl{\Vector{w}}$ denote the weighted lexicographic order given by the weight vector $\Vector{w}$. The monomial order $\wlexl{\Vector{w}}$ induces the following order on the polynomial variables: 
$$ X_1 \wlexl{\Vector{w}} X_2 \wlexl{\Vector{w}} X_3 \wlexl{\Vector{w}} X_4,$$
which is the same variable order induced by $\lexl$ in \Cref{Example:LexOrder}. Although the order of the polynomial variables is the same, the monomial order $\wlexl{\Vector{w}}$ is not the lexicographic order. For example, the four monomials given in \Cref{Example:LexOrder} is ordered as follows with respect to $\wlexl{\Vector{w}}$:
$$ X_2X_3 \wlexl{\Vector{w}} X_1 \wlexl{\Vector{w}} X_2X_4 \wlexl{\Vector{w}} X_4.$$
\end{example}

\begin{remark}
\label{Remark:MonOrder}
By changing the weight vector $\Vector{w}$ in \Cref{Example:WLexOrder}, one can obtain other monomial orders which induce a variable order different than the ones induced by the monomial orders given in \Cref{Example:LexOrder} and \Cref{Example:WLexOrder}. As an example, consider the new weight vector $\Vector{w'} = (1,2,3,4)$. The induced variable order in this case is
$$ X_4 \wlexl{\Vector{w'}} X_3 \wlexl{\Vector{w'}} X_2 \wlexl{\Vector{w'}} X_1,$$
\end{remark}

Let $f \in \BB({\Vector{X}})$ be a nonzero polynomial, $f_{\Vector{\alpha}}$ denote the coefficient of the monomial $\Vector{X}^{\Vector{\alpha}}$ in $f$ for $\Vector{\alpha} \in \mathbb{Z}^n_{\geq 0}$, and 
$$A_f:=\left\{ \Vector{\alpha}\in \mathbb{Z}^n_{\geq 0} \  \mid \ f_{\Vector{\alpha}} \neq 0  \right\}$$
denote \emph{the support of f}. 
If $>$ is a monomial order on $\BB({\Vector{X}})$, then \emph{the multidegree of f} is 
$$ \multideg{f} :=  \underset{\Vector{\alpha} \in A_f}{\arg\max} \left\{ \Vector{X}^{\Vector{\alpha} }\right\},$$
where the maximization is taken with respect to the monomial order $>$. Moreover, \emph{the leading monomial of $f$} is $\LM{f} := \Vector{X}^ {\multideg{f}}$, the \emph{the leading coefficient of $f$} is $\LC{f} := f_{\multideg{f} }$, and the \emph{the leading term of $f$} is $\LT{f} := \LC{f}\cdot \LM{f}$.
If $G$ is a finite set of polynomials in $\BB({\Vector{X}})$, then $\LM{G}$ denotes the set containing the leading monomials of the polynomials in $G$.

\begin{definition}
		Let $\cal I$ be an ideal of the Boolean polynomial ring $\BB(\Vector{X})$, a finite subset $G\subset \cal I$ is called a \emph{Boolean Gröbner basis} of $\cal I$ if $\left\langle \LM{G} \right\rangle = \left\langle \LM{\cal I} \right\rangle$.
\end{definition}

Recall that we represent each observation with a Boolean polynomial in $\BB(\Vector{X})$. We utilize Gröbner basis of an ideal $\cal I$ to simplify these polynomials, and hence simplify the rows encoded by the Boolean data that we started with.

\begin{proposition}[\cite{CoxLittleOshea:Ideals}]
    \label{Proposition:GroebnerMain}
    Let $\cal I$ be an ideal of $\BB(\Vector{X})$, $G = \{g_1, \dots, g_k \}$ be a Gröbner basis of $\cal I$, and $f$ be polynomial in $\BB(\Vector{X})$. Then there exists a unique $r \in \BB(\Vector{X})$ such that:
    \begin{enumerate}
        \item $f = g + r$ for some $g \in \cal I$,
        \item Non of the polynomials $\LT{g_1}, \dots, \LT{g_k}$ divide any term of $r$.
    \end{enumerate}
\end{proposition}

\begin{definition}
\label{Definition:NormalForm}
Given an ideal $\cal I$ of $\BB(\Vector{X})$, a Gröbner basis $G$ of $\cal I$, and a polynomial $f\in \BB(\Vector{X})$, then 
the polynomial $r$ given by the \Cref{Proposition:GroebnerMain} is called as the \emph{normal form of $f$ with respect to the Gröbner basis $G$}. Note that $G$, and consequently the normal form $r$, depend on the choice of monomial ordering on $\BB(\Vector{X})$.
\end{definition} 

One can think of the normal form $r$ as a reduction of the original polynomial $f$ obtained by setting each polynomial in the ideal $\cal I$ equal to zero. In order to compute the $r$, one needs to perform a series of polynomial divisions with remainder; first divide $f$ by $g_1$ to obtain the remainder $r_1$, and then at each next step divide the remainder $r_{i-1}$ by $g_{i}$ successively. 
\begin{remark}
\label{Remark:GroebnerBasisUniqueNormal}
If the set $ G = \{g_1,\dots,g_k\}$ is not a Gröbner basis of the ideal $\cal I$, then the result of taking successive remainders depend on the order of polynomials in $G$. However, if $G$ is a Gröbner basis of the ideal that it generates, then the result do not depend on $g_i$'s order. 
\end{remark}

Now, we algorithmically describe the process we apply in order to obtain a reduction scheme from a given Boolean data table. We start with a Boolean data table given as in \Cref{Definition:BooleanDataTable}, and we associate the Boolean variables $X_1,\dots,X_n$ to each of its columns from left to right, respectively. In order to compute normal forms we have to fix a monomial order as discussed in \Cref{Remark:GroebnerBasisUniqueNormal}. We exclusively use the weighted lexicographical order defined in \Cref{Definition:WLexOrder}, where the weight of the variables are computed as in \Cref{Algorithm:Weight}. This algorithm computes the weights of each variable $X_i$ by considering the relation between occurrences of $true$ and $false$ for this variable in the Boolean data table $M$. The weight is 
$$
w[i]=n_{true}\cdot n_{false}.
$$
In our perspective, $X_i$ is more likely to be a distinctive property either if it occurs in only some of the observations, or if it occurs in most of the observations. In this case the weight in Algorithm 1 would be small.

In terms of the algebraic reduction, variables $X_i$ with higher weight are less likely to appear in the normal forms we compute than the variables with lower weights. In terms of the data, a high $w[i]$ shows that nearly half of our observations admit the characteristic $X_i$, and indicates that the characteristic $X_i$ is not biased towards a minority in the set of observations. Therefore, the variables that correspond to a less biased characteristics are less likely to appear in the reduced expression obtained by taking the normal form. 
\begin{remark}
The weight vector computed with \Cref{Algorithm:Weight} is a modular component of our main reduction described in \Cref{sec:simplify}. Therefore, the ideas we present in the next section can be replicated for a different weight assignment scheme.
\end{remark}

\begin{algorithm}[H]
\caption{An algorithm for determining variable weights}\label{Algorithm:Weight}
\begin{algorithmic}[1]
\Require{Boolean data table: $M \in \BB^{ m \times n}$} 
\Ensure{A vector of weights:  $n\in \mathbb{Z}^n$ }
\State $i \gets 1$
\While{$i \leq n$}
    \State $n_{true} \gets 0$
    \State $j \gets 1$
    \While{$ j \leq m$}
        \If{$M[j,i]\, is \,  true $}
            \State $n_{true} \gets n_{true}+1$
        \EndIf
        \State $j \gets j+1$
    \EndWhile
    \State $w[i] \gets n_{true}\cdot(m-n_{true})$
    \State $i \gets i+1$
\EndWhile
\State \Return $w$
\end{algorithmic}
\end{algorithm}

\begin{remark}
In fact, it is even possible to other monomial orderings instead of the weighted lexicographical order we propose to use. There are various ways to define a monomial order on $\BB(X_1,\dots,X_n)$ other than the ones we introduced here, see e.g., \cite{CoxLittleOshea:Ideals}. 
\end{remark}

\subsection{Finding Peculiarities} \label{sec:simplify}
The problem which motivated this article is not finding logical rules by trusting the data set, but to find possible inconsistencies in the data set. Imagine there are ``simple'' (yet unknown) logical rules represented by the reported observations (of a large data set). We describe a process in \Cref{Algorithm:Reduction} to find an alternative reduced formulation of the polynomials corresponding to the objects. 

We note that the method given in \Cref{Algorithm:Reduction} can be implemented in many modern computer algebra systems. We require two key functionalities from the computer algebra system, which are to compute a Gröbner basis and to compute the normal form with respect to a Gröbner basis. 
\begin{remark}
\label{Remark:Website}
In our explicit computations, we use the Gröbner basis implementation available in the mathematical software OSCAR \cite{OSCAR}, and our code is available in the following link:
\begin{center}
\small
    \texttt{https://github.com/OguzhanYueruek/AlgebraicClusteringSupplementaryData}
\end{center}
\end{remark}

The crucial step in our analysis of a given binary table $M$ is to generate the Aclus ideal $\cal I$ associated to $M$. 
In order to compute a generator for $\cal I$, we first compute the Boolean polynomials $g_1,\ldots,g_m$ corresponding the rows of $M$, and consider the Boolean polynomial $g := 1 + \sum_{i=1}^m g_i $. Then, the Aclus ideal associated to $M$ is generated by the polynomial $g \in \BB(\Vector{X})$. 

\begin{algorithm}[H]
\caption{Algebraic Clustering ({\sc Aclus}) algorithm for reducing an observation with respect to a given Boolean data table}\label{Algorithm:Reduction}
\begin{algorithmic}[1]
\Require{Boolean data table: $M \in \BB^{ m \times n}$, Vector of weights:  $\Vector{w}\in \mathbb{Z}^n$, Boolean polynomial to reduce: $f \in \BB(X_1,\dots,X_n)$} 
\Ensure{A reduced polynomial $r \in \BB(X_1,\dots,X_n)$}
\State $g \gets 1$
\State $j \gets 1$
\While{$j \leq m$}
    \State $i \gets 1$
    \State $expression \gets 1$
    \While{$ i \leq n$}
        \If{$M[j,i] \, is \, true $}
            \State $expression \gets expression \cdot X_i$
        \Else
            \State $expression \gets expression \cdot (1+X_i)$
        \EndIf
        \State $i \gets i+1$
    \EndWhile
    \State $ g \gets g + expression$
    \State $j \gets j+1$
\EndWhile
\State $I \gets \left\langle g \right\rangle  $ \Comment{Aclus ideal generated by $g$ in $\BB(X_1,\dots,X_n)$} 
\State $G \gets groebner\_basis(I,wlex,w)$ \Comment{ {\scriptsize Gr\"obner basis w.r.t the wlex order with weight $w$.}}
\State $r \gets normal\_form(f,G)$ \Comment{Normal form of $f$ w.r.t. $G$, see \Cref{Definition:NormalForm}} 
\State \Return $r$
\end{algorithmic}
\end{algorithm}

Note that taking the union of all observation patterns that satisfy $g_i$ and $g_j$, in Boolean algebraic terms, is equal to the set of all observation patterns that satisfy $g_i + g_j + g_ig_j$. However, since the rules given by each $g_i$ is an atom in the Boolean ring $\BB[\Vector{X}]$, for any $1\leq i,j \leq m$ we have $g_ig_j = 0$. Therefore, in order to compute the union all observation patterns that satisfy every $g_i$, it is enough to consider the rule given by $g_1 + \dots + g_m$.  By considering the polynomial $g = 1 + g_1 \dots + g_m$, we find a rule that is not satisfied by any of the observation patterns given by the rows of the binary table $M$.

What happens, if a coder has added a wrong observation to this binary data set? In this case,  algebraic expressions which should be part of the Aclus ideal $\cal I$ are taken away from $\cal I$, because such an expression is now part of the table. This wrong observation cancelled a rule which in reality applies. In order to make this rule visible, we take the Boolean polynomial, which represents this observation and compute its remainder with respect to the ideal, i.e. its representation in the quotient ring $\quot{\BB(\Vector{X})}{\cal I}$.

\begin{example}
The outcome of applying the reduction method to the polynomials which represent the lines of Table \ref{tab:my_label} is given in \eqref{eqn:my_label} for its 18 rows. The details of the computation can be found in the jupyter notebook called \texttt{AlgebraicClustering\_Tree\_weighted.ipynb}, see \Cref{Remark:Website}.
\begin{align}\label{eqn:my_label}
\nonumber
1 & \rightarrow  w a g + w g + a c + a + c g + g \\ \nonumber
2 & \rightarrow  a c + a g + a + n v g + n c + n g + n + v g + c + g + 1 \\ \nonumber
3 & \rightarrow  n v g + v g \\ \nonumber
4 & \rightarrow  c g + c + g + 1 \\ \nonumber
5 & \rightarrow  n v g + n c + n g + n + c g + c + g + 1 \\ \nonumber
6 & \rightarrow  n v g \\ \nonumber
7 & \rightarrow  w a g + a c + a \\ \nonumber
8 & \rightarrow  a c + a g + a \\ \nonumber
9 & \rightarrow  w a g + w a + a g + a \\ 
10 & \rightarrow  w a g + w g + n c + n + c g + c + g + 1 \\ \nonumber
11 & \rightarrow  a v + a g + a + n v g + n v + n g + n + v g + v + g + 1 \\ \nonumber
12 & \rightarrow  w a g + w a + w g + w + a v + a g + a + n v g + n v + v g + v + c g + c \\ \nonumber
13 & \rightarrow  n c + n + c g + c + g + 1 \\ \nonumber
14 & \rightarrow  n v g + n v + n g + n + c g + c + g + 1 \\ \nonumber
15 & \rightarrow  n v g + n v \\ \nonumber
16 & \rightarrow  w a g \\ \nonumber
17 & \rightarrow  a v + a g + a \\ \nonumber
18 & \rightarrow  w a g + w a + a v + a g + a \nonumber
\end{align}

Here, the expression $cg+c+g+1$ discussed in Section \ref{sec:intro} becomes visible for line No. 4 of the data set. Note that ``$cg+c+g$'' provides the union set of all observations with property $c$ and of all observations with property $g$. The expression ``$\ldots+1$'' denotes the complement of this union set.  Therefore, this expression can also be formalized as $\lnot(c\vee g)$.

Another example: In line No. 6 the expression ``$nvg$'' denotes the observation that at a very cold weather condition ($v$) in a neutral soil ($n$) the tree grows well. This might be a conspicuous observation, because in line No. 15 the tree does not grow in neutral soil at a very cold weather condition. In the same sense line No. 3 ``$nvg+vg$'' might be a conspicuous observation, because here $nvg+vg$ can be transformed into $(n+1)vg$ which means that in a non-neutral soil at very cold weather condition the tree grows well.   

Some expressions are more complicated. The logical interpretation of these expressions needs some practise. The meaning of expression ``$ac+ag+a$'' in line No. 8 becomes more clear, if we transform it into $a(c+g+1)$. Being a product it means that the relevant conditions occur together. Thus, line No. 8 is an observation in acidic environment ``$a$'' {\em and} at the condition ``$1+c+g$''. The latter one is only valid if {\em either} $c$ and $g$ both hold {\em or} both $c$ and $g$ do not hold, which can be formalized as $c\Leftrightarrow g$. Interpretation:  The observation in line No. 8 is special, because at acidic conditions the setting of $c$ is equal to the setting of $g$. This is not true for any other row of Table \ref{tab:my_label}. If line No. 8 would not be part of the observations, then at acidic conditions, $g$ and $c$ would never have the same value (the copper beech either grows well or it is cold).  A similar algebraic transform and a similar interpretation is possible for line No. 17.
 
Another example of how to interpret the algebraic expressions is given by line No. 9. Here, we can transform the expression ``$wag+wa+ag+a$'' into ``$a(g+1)(w+1)$''. Again it is a product. It is a statement about acidic soil. ``$g+1$'' means that it is a situation in which the copper beech does not grow well. In this situation line No. 9 has the special feature ``$w+1$'', which is true if $w$ does not hold. Interpretation:  The observation No. 9 is special, because the copper beech does not grow well in acidic soil and the soil is not wet. This is not true for all other rows (with $\lnot{g}$ and with $a$) in Table \ref{tab:my_label}. If observation No. 9 would not be part of the data set, then ``not growing well in acidic soil'' means that it would be wet. 
\end{example}

Some algebraic expressions are more complicated to understand, because they (as a single row) do not contradict ``simple logical rules''. This might be a good indicator that in these rows the characteristics are correctly assigned. However, how can a coder see whether one has to combine some rows to form a cluster of lines with possibly wrong entries? 

\begin{example}
Have a look at the simplified expressions in Equation \eqref{eqn:my_label}. The expression in line No. 3 is ``$nvg+vg$'', whereas line No. 6 is ``$nvg$''. Adding these two expressions with the rules ($x+x=0$) of Boolean rings leads to cancellation of $nvg$. The remaining expression is $vg$. This expression denotes the peculiarity of the selected cluster $\{3,6\}$ of rows of the table which has already been discussed in Section \ref{sec:intro}: $v$ and $g$ occur together only in these rows. It is possible to represent the union of these two rows by simply summing up the polynomial expressions, because the product of $nvg+vg$ and $nvg$ is equivalent to zero with regard to the quotient ring: The different expressions in \eqref{eqn:my_label} are the result of a mapping from the ``full'' polynomials (representing the observations) to the quotient ring. The product of two ``full'' polynomials is zero, because two different patterns at least differ in one aspect. Thus, $(nvg+vg)nvg=0$ and therefore $(nvg+vg) \cup nvg = nvg+vg+nvg=vg$.

It is not easy to find other combinations of two rows such that adding their terms leads to a ``simple'' logic expression including the variable $g$. Another example for such a combination are line No. 7 and line No. 9. Adding their algebraic expressions leads to: $a(w+c+g)$. If one would leave out these two lines in the table, then $a(w+c+g)=0$. In this case, an acidic environment means that exactly two of the three variables $w$, $c$, and $g$ are always set (have a look at the table). The additional rule for an acidic soil which then would apply can also be read from the algebraic expression $w+c+g=0$: If it is cold and wet, the beech does not grow. If it is either cold or wet then the beech grows. 

Combining two rows only exceptionally leads to ``simple'' algebraic expressions including $g$. This is an indicator for the assumption that other clusters of experiments are correctly reported. They do not contradict ``simple'' rules. 
\end{example}

\section{Real-World Examples}

We have used the approach of this article to examine given binary data tables (up to 1500 rows and up to 19 columns) for possible inconsistencies.

\subsection{Questionnaire}
As a first example, we took the raw data from a questionnaire about housing decisions of senior Australians \cite{Aust}. 

One question to the 1524 participants was ``Which of the following were considerations that influenced when you retired?'' and possible answers were:
\begin{itemize}
    \item[e] Having (e)nough money to live comfortably for the full length of your retirement
   \item[h] Having to sell your (h)ome in order to fund your retirement
   \item[I] Being able to fund your retirement and still leave an (I)nheritance to your children
   \item[a] Your (a)ge and its impact on your ability to access the age
   \item[c] Your assets and in(c)ome and their impact on your 
   age
   \item[p] Having enough money to afford your desired housing arrangements e.g. a home of your choice, your (p)referred retirement village etc.
   \item[H] Your (H)ealth or the health of your partner
   \item[d] Your ability to (d)raw down upon your superannuation
    \item[x]Job loss/Loss of employment
\end{itemize}
This has led to a $1524\times 9$ binary table, which has been analysed with our algebraic approach {\sc Aclus}. The specific computation done with this data set can be found in the file \texttt{AlgebraicClustering\_Retirement\_weighted.ipynb}, see \Cref{Remark:Website}. After computing the ideal $\cal I$ on the basis of the $9$ listed Boolean variables, we have reduced the information included in the answers of the $1524$ individuals. Some answer patterns then revealed to be remarkable. 

Interesting is the answer pattern of person No. 1104. This pattern reduced to $hacd(p+1)$ with respect to the Aclus ideal $\cal I$. This person had to sell home in order to fund retirement, answered ``yes'' to the questions ($a$ and $c$) about the impact on age, and this person can draw down upon superannuation. Interestingly, this person did not answer ``yes'' to ``Having enough money to afford your desired housing arrangements e.g. a home of your choice, your preferred retirement village'', which could mean that selling home was not enough to can take these aspects into account.

Two individuals (No. 913 and 1371) have a reduced pattern $Iapdx(c+1)$. Reading the corresponding answers, it is really astonishing that ``Your assets and income and their impact on your age'' has not been answered with ``Yes'', because financial aspects seemed important for these two persons. 

If we combine the pattern $Iacpx + Iacdx + Icpx + IcHdx$ of person No. 866 with the pattern $Iacpx + Iacdx$ from person No. 875, then the sum is $Icx(p+Hd)$. The two persons have lost their jobs, they still want to have enough to inherit to their children, and important for them are their assets and income. The peculiarity is (reading that they answered yes to $p$ and $H$) that they not answered ``yes'' to ``Your ability to draw down upon your superannuation''. The loss of their jobs seemed to have been a quite severe event. 

The reduction of the answer patterns of each individual person to the ``remarkable'' part can be used to analyse questionnaires in order to eventually dig deeper into the situation of the person who gave the answers. Inherent rules of the answers can be read from the Gr\"obner basis of the ideal $\cal I$. In the case of this questionnaire, one basis element is $hIcx + hcx \in {\cal I}$. This means $0=hIcx + hcx$, or equivalently $hcx = hIcx$. Every person who had to sell home, who said that assets and income have been important, and who lost job also answered that being able to fund retirement and still leave an inheritance to the children has been an influencing consideration.

\subsection{Egyptian Statues}
Another sample datasheet contained classifying properties of a set of ancient Egyptian statues which were found in the so called Cachette of Karnak (see \cite{Coulon2016}). It was provided by Egyptologist Ralph Birk (Freie Universität Berlin), compiled with help of Sarah M. Klasse. For each statue, it was reported whether or not certain characteristics applied. In our example, metadata on 495 statues has been extracted from the database  “Cachette de Karnak” \cite{karnak}, including title of the object, excavation dates, statue type, material and measurements. This set has been then enriched with further categories (e.g. who is depicted, a god or a private person?) and an inductively constructed classification of destruction types (e.g. is the head amputated? See \cite{Jambon2016}). Given the vast array of individal, highly specific objects, our colleagues were first interested in testing their intracoder reliability, e.g. the consistency in how the same person coded the data through time \cite{Connor2020}. Second, they were interested in refining their coding and their research questions in respect to the chosen categories after a review of the generated reductions of {\sc Aclus}.

In the binary table in our supplemental material, see \Cref{Remark:Website}, we have deleted basic information about the statues and many of the properties in order to anonymize the original data set. Following 6 characteristics remained for illustration reasons:
\begin{itemize}
    \item[A] Amputations
   \item[h] amputated head
   \item[H] missing head
   \item[a] partial amputation
   \item[x] partially destroyed
   \item[X] largely destroyed
\end{itemize}

In total 29 different classification patterns have been found among the 495 statues.
From taking the complement of the sum of the polynomials corresponding to these patterns the ideal ${\cal I}$ has been generated. This is the ideal of ``rules'' which govern the dataset of statues. The most frequent pattern was “none of the destruction properties applies” (174 times). The Gr\"obner basis of the ideal only has (besides the trivial expression like $A^2 + A$) ten generating ``rules'':
\begin{eqnarray*}
 && hX + HX, \cr
 && AHX + HX, \cr
 && HaX + HxX, \cr
 && hHx + Hx, \cr
 && AxX + xX, \cr
 && AhH + AH + hH + H,\cr
 && hHa + Ha, \cr
 && AHx + Hx, \cr
 && AHa + AH + hHa + H, \cr
 && AhHa + AhHx + AhH + Ahax + Ahx + hHa + hHx + hH + hax + hx.
\end{eqnarray*}
Some of these rules are human understandable.  One generating element is $hX+HX\in{\cal I}$. $hX+HX=0$ means that $hX=HX$, i.e., a statue which is largely destroyed and has an amputated head, also has its head missing. By $AHX+HX$ we see that these statues also have the property ``amputations''. $hHx+Hx=0$ means that a statue which has a missing head and is partially destroyed also has an amputated head. By $AHx+Hx$ it further follows that these statues have ``amputations''. 

Now we took this ideal and reduced the 29 different destruction patterns according to it. In this case some reduced patterns are remarkable. One reduction is $(A+1)H$. This is one statue which does not have an $A$mputation but a missing  $H$ead. Indeed this had been a false classification detected by {\sc Aclus}. Another polynomial has been reduced to $(h+1)H$, a statue with a missing head ($H$) which does not have an amputated head ($h$). This classification was also revised after detection. 

Taking all properties into account, another example was $Tb+T$: a statue which is not a $b$lock statue but has an intended destruction of $T$ext. This was not a mistake but a peculiarity of that statue. However, inspired by the unexpected rule ($Tb+T=0$ would mean that $b\succeq T$: all intended destructions of texts belong to block statues) we had a closer look to the block statues and found one statue with destructed text which had been classified as ``not intended''. After deeper analysis of the surface of the statue, the classification had been changed to ``intended destruction of text''.  

\paragraph{Preliminary results and refinement of classifications.} The given dataset is only a subset in an ongoing research project. We analysed here, as a second example, the information on the destruction patterns and the social status of the depicted person (private, royal, divine). In total 19 different binary categories had been investigated for the individual statues. Some binary patterns occurred more often, such that 226 different patterns remained. From these patterns we computed the ideal $\cal I$. Computing the Gr\"obner basis of this ideal revealed some rules. One rule, e.g. is $pg=0$, which was the code for $p$rivate and $g$od. This rule means that the categories ``depicting a private person'' and ``being a divine statue'' are exclusive. This is a rule which our colleagues implemented in the classification scheme, but was found purely by analysing the binary table with algebraic methods, namely by stating that $pg\in {\cal I}$. Another rule was $ck\in{\cal I}$, which meant that there is no $c$omplete royal ($k$ing) statue in the given selection of 495 statues, as well as $HaGb\in{\cal I}$: There is no granite block statue with partial amputations together with a missing head.  

Both point to systematic qualities of the studied set which need further investigation.

\section{Conclusion}
Our method {\sc Aclus} helped researchers from Egyptology to refine their coding scheme for classifying ancient statues. The proposed method extracts all possible logical rules about the properties of observed objects in terms of an ideal of Boolean Polynomials. The Gr\"obner Basis of this ideal depends on the monomial ordering. In this article, we decided for using monomial orderings based on the number of occurrences of characteristics. One particular order called the \emph{elimination order} would be another good candidate which offers various functionalities, which we plan to discuss in a follow up work.

By using the rules given by the ideal, we can reduce the information, i.e., the polynomial which represents the characteristics of one object -- included in one row of the binary data table -- to a shorter algebraic expression. This expression provides the peculiarity of the corresponding object in terms of a new rule which would apply if that object would be deleted from the data table. 

Only as a remark about the illustrative example of this article: Table \ref{tab:my_label} has been created to represent the properties of the copper beech according to the decision tree in \cite{rotbuche}. However, it is obvious that this decision tree is not correct. The author based decisions on a difference between the German word ``basisch'' and the German word ``alkalisch'' which actually means the same (namely high pH). Thus, we can expect that corrections in Table \ref{tab:my_label} lead to more simple rules about copper beeches (like in ``Rotbuchen'', \cite{ullrich}). Maybe one rule about copper beeches is even as simple as $\overline{c}\Rightarrow g$ (``If it is not cold, then the copper beech grows well.''). In this case only line No. 4 would be a wrongly reported observation found by {\sc Aclus} and would have to be deleted from the table in order to correct it. 

\paragraph{Conceptual aspects.} Based on the algorithm {\sc Aclus}, the consequences of removing and adding observations to an existing binary data set can be distinguished. There are the following situations:
\begin{itemize}
    \item[1a] {\em Adding an observation with a new pattern.} In this case, elements of $\cal I$ are removed. There are now less possible logical rules about the relations between the properties of the objects (in terms of a subset relation). 
    \item[1b] {\em Adding an observation with a pattern that already exists.} In this case, the logical rules do not change, but the importance of certain properties can change their weighting. Therefore, the representation of rules and residues are possibly changing due to a possibly new Gr\"obner basis.  
    \item[2a] {\em Removing an observation with a unique pattern.} In this case, the removed pattern is added to $\cal I$. There are now more possible logical rules about the relations between the properties of the objects (superset relation). 
    \item[2b] {\em Removing an observation with a pattern that is repeated.} Same as 1b.
\end{itemize}

\subsection*{Acknowledgements}  We thank Ralph Birk for providing information about ancient Egyptian statues and for fruitful discussions about our computational results. We want to thank Konstantin Fackeldey and Michael Joswig for fruitful conceptual discussions with regard to {\sc Aclus}. The work from M.W. is part of the MATH+ project ``EF5-4 -- The Evolution of Ancient Egyptian – Quantitative and Non-Quantitative Mathematical Linguistics''. O.Y. is also funded by the Deutsche Forsch\-ungsgemeinschaft (DFG, German Research Foundation)  under Germany's Excellence Strategy - The Berlin Mathematics Research Center MATH+ (EXC-2046/1, project ID 390685689, sub-project AA1-9).

\printbibliography

\end{document}